\documentclass{SLC}  
\articlenumber{2}

\newcommand{\BE}{\begin{equation}}	
\newcommand{\EE}{\end{equation}}
\def\OEIS{{\footnotesize OEIS}}

\author[A.~Guttmann]{Anthony J. Guttmann\textsuperscript{1}\protect\orcid{0000-0003-2209-7192}}
\address{\textsuperscript{1}School of Mathematics and Statistics, The University of Melbourne,
Vic.~3010, Australia;  \website{https://blogs.unimelb.edu.au/tony-guttmann}}
\author[V.~Kot\v{e}\v{s}ovec]{V\'aclav Kot\v{e}\v{s}ovec\textsuperscript{2}\protect\orcid{0000-0003-4516-0617}}
\address{\textsuperscript{2}Prague, Czech Republic; \website{http://www.kotesovec.cz/math.htm}}

\title{A numerical study of $L$-convex polyominoes and 201-avoiding ascent sequences}

\abstract
{For $L$-convex polyominoes we give the conjectured asymptotics of the generating function coefficients, obtained by analysis of the coefficients derived from the functional equation given by Castiglione et al. 

For 201-avoiding ascent sequences, we conjecture the solution, obtained from the first twenty-three coefficients of the generating function. This solution is D-finite, indeed algebraic. The conjectured
solution then correctly generates all subsequent coefficients. We also obtain the asymptotics, both from direct analysis of the coefficients, and from the conjectured solution.

As well as presenting these new results, our purpose is to illustrate the methods used, so that they may be more widely applied.}

 
\begin{document}
\maketitle

\section{Introduction}
\label{introduction}
\begin{figure}[ht]
 \centerline{  \includegraphics[width=0.35\linewidth]{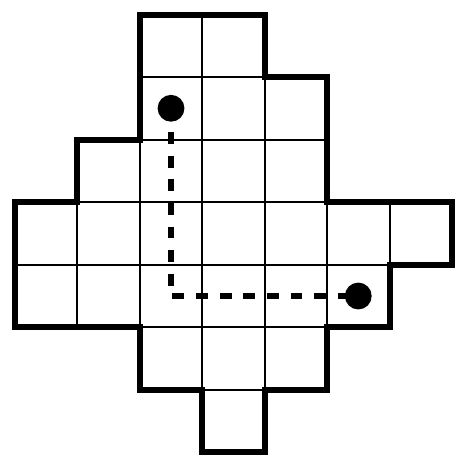} }
   \caption{An $L$-convex polyomino.}
 \label{Lconvex}
\end{figure}

In~\cite{CFMRR7}, Castiglione et al.~gave a functional equation for the number of $L$-convex polyominoes. These are defined as polyominoes with the property that any two cells may be joined by an $L$-shaped path, that is to say, a path with at most one right-angle bend. An example is shown in Figure~\ref{Lconvex}.
It can be seen that such polygons can be described as a stack polyomino placed atop an upside-down stack polyomino. A stack polyomino is just a row-convex bargraph polyomino. The perimeter generating function of $L$-convex polyominoes has a simple, rational expression, 
\BE P(x)=\frac{(1-x)^2}{2(1-x)^2-1} = 1+2x+7x^2+24x^3+\cdots,\EE
and is the sequence~\oeis{A003480} in the On-line Encyclopaedia of Integer Sequences (OEIS),~\cite{OEIS}.
Accordingly, one has \BE[x^n]P(x) =\frac{(2+\sqrt{2})^{n+1}-(2-\sqrt{2})^{n+1}}{4\sqrt{2}} \sim \frac{1+\sqrt{2}}{4}(2+\sqrt{2})^n.\EE
The area generating function is given by Castiglione et al.~\cite{CFMRR7}  
\BE \label{eq:lconv}
A(q)=1+\sum_{k \ge 0} \frac{q^{k+1}f_k(q)}{(1-q)^2(1-q^2)^2\cdots(1-q^k)^2(1-q^{k+1})} = 1+q+2q^2+6q^3+15q^4+ \cdots,
\EE
 where \BE f_k(q)=2f_{k-1}(q)-(1-q^k)^2 f_{k-2},\EE
with initial conditions $f_0(q)=1$, and $f_1(q)=1+2q-q^2$. 
We used this expression to generate 2000 terms of the sequence, and these are given in the OEIS as sequence \oeis{A126764}. Analysis of this sequence allowed us to derive the conjectured asymptotics as
\BE [q^n]A(q) \sim \frac{13\sqrt{2}}{768\cdot n^{3/2}}\exp(\pi\sqrt{13n/6}).\EE 
 In the next section we will describe how this estimate was obtained.

The second problem we are considering is that of 201-avoiding ascent sequences, defined below.
  Given a sequence of non-negative integers, $n_1 n_2 n_3 \ldots n_k$,  the number of {\em ascents} in this sequence is \BE asc(n_1 n_2 n_3 \ldots n_k) = |\{ 1  \le j<i : n_j \le n_{j+1} \} |.\EE
  
  The given sequence is an {\em ascent sequence} of length $k$ if it satisfies $n_1=0$ and 
\BE n_i \in [0,1+asc(n_1 n_2 n_3 \ldots n_{k-1} )]  \text{ for all } 2 \le i \le k.\EE 
For example, $(0,1,0,2,3,1,0,2)$ is an ascent sequence, 
  but $(0,1,2,1,4,3)$ is not, as $4 > asc(0121)+1=3$. 
  
  Ascent sequences came to prominence when Bousquet-M\'elou et al.~\cite{BCDK10} related them to $(2+2)$-free posets. They have subsequently been linked to other combinatorial structures. See~\cite{K11} for a number of examples. Later, Duncan and Steingr{\'{\i}}msson~\cite{DS11} studied {\em pattern-avoiding ascent sequences. }
  \pagebreak

  A pattern is simply a word on nonegative integers (repetitions allowed). Given an ascent sequence $(n_1 n_2 n_3 \ldots n_k)$,  a pattern $p$ is a subsequence $n_{i_1}n_{i_2}\ldots n_{i_j},$ where $j$ is just the length of $p$, and where the letters appear in the same relative order of size as those in $p$. For example, the ascent sequence $(0,1,0,2,3,1)$ has three occurrences of the sequence $001,$ namely $002,$ $003$ and $001$. If an ascent sequence does not contain a given pattern, it is said to be {\em pattern avoiding}. 
 \smallskip
 
 The connection between pattern-avoiding ascent sequences and other combinatorial objects, such as set partitions, is the subject of~\cite{DS11}, 
while the connection between pattern-avoiding ascent sequences and a number of stack sorting problems is explored in~\cite{CCF20}.
 
 Considering patterns of length three, the number of ascent sequences of length $n$ avoiding the patterns $001$, $010$, $011$, and $012$ is $2^{n-1}$ (the sequence \oeis{A000079} in the OEIS). 
For the pattern $102$ the number is $(3^n+1)/2$ (\OEIS\ \oeis{A007051}), while for $101$ and $021$ the number is just given by the $n^{th}$ Catalan number (\OEIS\ \oeis{A000108}).
\smallskip
  
For the pattern 201, the first twenty-eight terms of the generating function are given in the OEIS as sequence \oeis{A202062}, and it is this sequence that we have used in our investigation.  
 First, we found, experimentally, that the coefficients given in the OEIS satisfied a recurrence relation, given in Section~\ref{201} below. 
 This recurrence can be converted to a second-order inhomogeneous ODE, or, as we prefer, a third-order homogeneous ODE.  The smallest root of the polynomial multiplying the third derivative in the ODE is $x=0.1370633395\ldots$ and is the radius of convergence of the generating function, and of course the reciprocal of the growth constant $\mu=7.295896946 \ldots$
 \smallskip

  This ODE, readily converted into differential operator form, can be factored into the direct sum of two differential operators, one of first order and one of second order. The solution of the first order ODE is a rational function
 while the solution of the second turns out to satisfy a cubic algebraic equation.
 This can be solved by one's favourite computer algebra package (we give the solution below), and expanding this solution, and adding it to the expansion of the solution of the first-order ODE, gives the required expansion.
\smallskip
 
 This analysis required only the first 24 terms given in OEIS, so the correct prediction of the next four terms gives us confidence that this is indeed the exact solution. 
 Expanding this solution and analysing the coefficients, as described in Section~\ref{201},  leads to us conjecturing
 the following asymptotic behaviour for these coefficients:
 \begin{equation}
 u(n) \sim C\frac{\mu^n}{n^{9/2}},
 \end{equation}
 where  \BE \mu=\frac{14}{3}\cos\left (  \frac{\arccos(\frac{13}{14})} {3} \right ) + \frac{8}{3}\EE
 and \BE C=\frac{35}{16}\left (\frac{4107}{\pi} - \frac{84}{\pi}\sqrt{9289}\cos\left (\frac{\pi}{3} +\frac{1}{3} \arccos\left [\frac{255709\sqrt{9289}}{24653006} \right ]\right )\right )^{1/2}.\EE

In the next two sections we give the derivation of the results given above.
\pagebreak

\section{\texorpdfstring{$L$}{L}-convex polyominoes} \label{Lconv}
As mentioned above, a typical $L$-convex polyomino can be considered as a stack polyomino placed atop an upside-down stack polyomino. Stack polyominoes counted by area have generating function
\begin{equation}S(q)=\sum s_n q^n =\sum_{n \ge 1} \frac{q^n}{(q)_{n-1}(q)_n},\end{equation} where $(q)_n :=  \prod_{k=1}^n (1-q^k)$, and, as first shown by Auluck~\cite{A51}, one has
 \BE s_n \sim \frac{\exp(2\pi\sqrt{n/3})}{8\cdot 3^{3/4} \cdot n^{5/4}}. \EE
 
 Thus putting two such objects together, one would expect a similar expression for the asymptotic form of the coefficients of 
the generating function~\eqref{eq:lconv}, that is to say, an expression of the form \BE  \label{str_exp} l_n \sim \frac{\exp(a\pi n^\beta)}{c n^\delta},\EE  where
 we write $L(x)=\sum l_n x^n$ for the ordinary generating function of $L$-convex polyominoes. We expect both exponents $\beta$ and $\delta$ to be simple rationals, as for stack polyominoes, and the constants $a$ and $c$ to be products of integers and small fractional powers.
 
 The analysis of series with asymptotics of this type is described in detail in~\cite{G15} and we will not repeat that discussion here, but simply apply the methods described there.

\begin{figure}[b!] 
\centerline{\includegraphics[width=.35\textwidth]{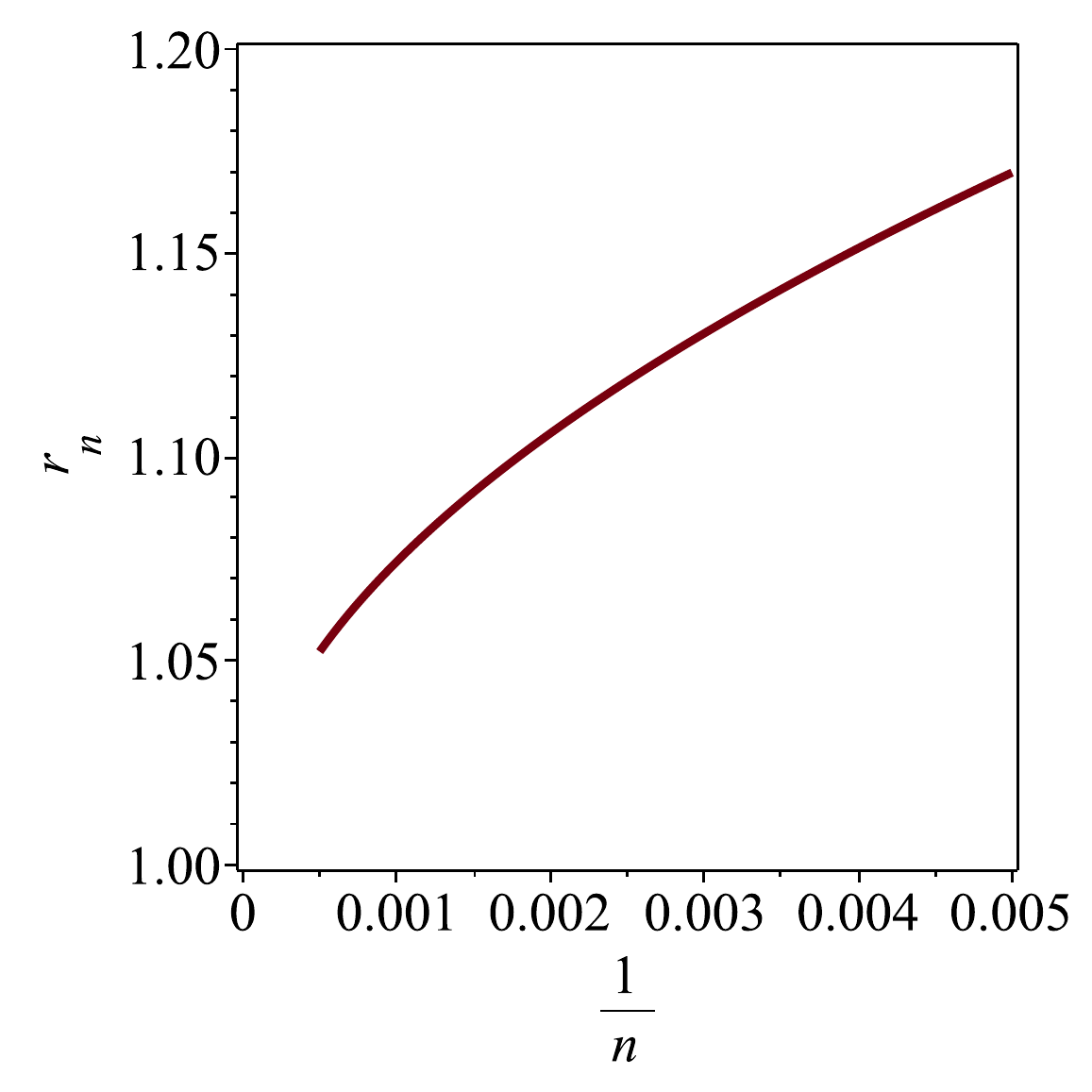}
\includegraphics[width=.35\textwidth]{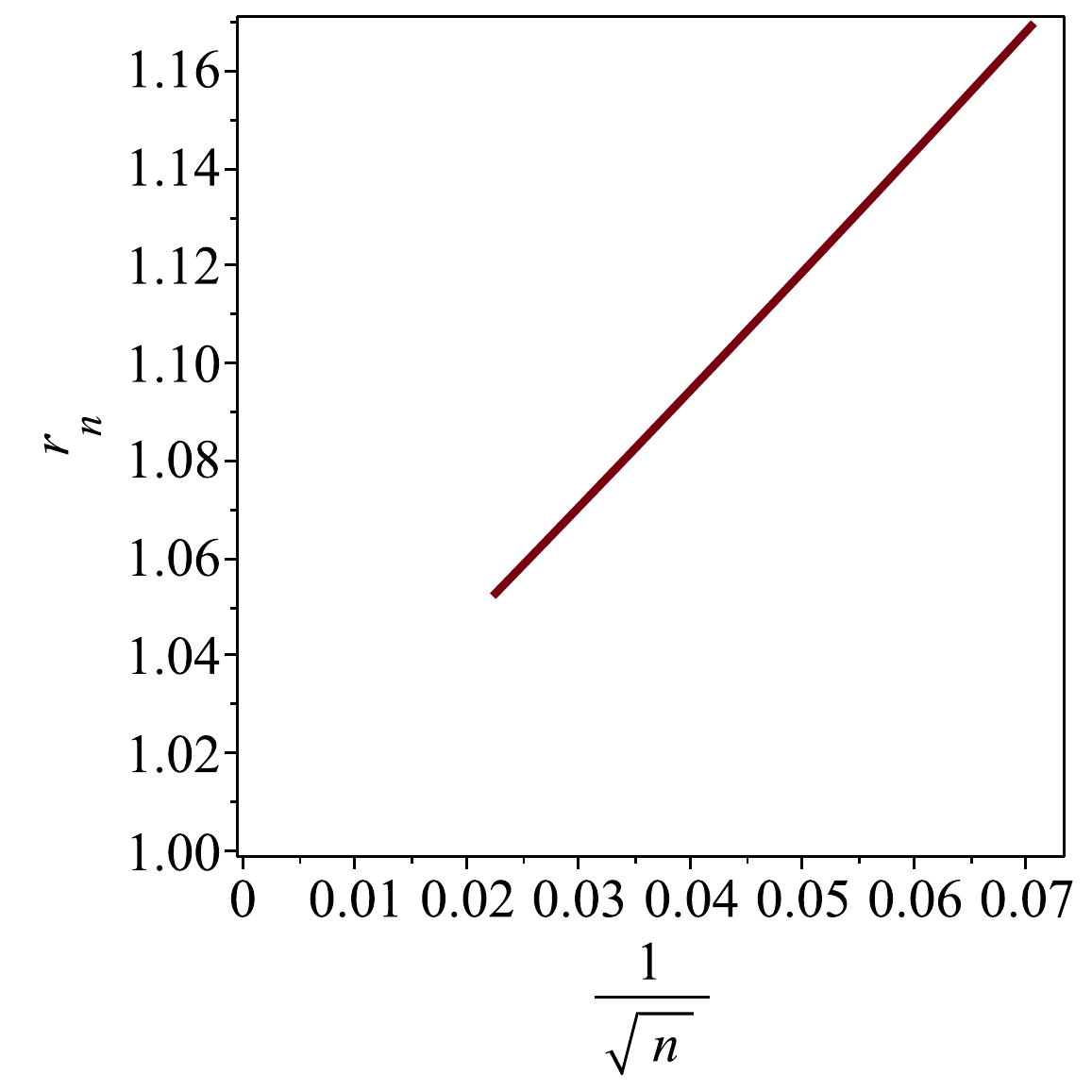}}\setlength{\abovecaptionskip}{0pt} 
\setlength{\belowcaptionskip}{0pt} 
\caption{$L$-convex ratios $r_n$ plotted against $1/n$ (left)  and against $1/\sqrt{n}$ (right).}
\label{fig:r12}
\end{figure}

 First, we consider the ratios of successive coefficients, $r_n = l_n/l_{n-1}$. For a power-law singularity, one expects the sequence of ratios to approach the growth constant linearly when plotted against $1/n$. In our case the growth constant is 1. That is to say, there is no exponential growth.
 From the asymptotic behaviour~\eqref{str_exp}, which is called of {\em stretched exponential} type,  it follows that the ratio of coefficients behaves as 
 \BE \label{eq:rn}
 r_n = \frac{l_n}{l_{n-1}} = 1+\frac{a\beta\pi}{n^{1-\beta}}+ O \left ( \frac{1}{n} \right ),
 \EE
  so we expect the ratios to approach a limit of 1 linearly when plotted against $1/n^{1-\beta}$, and to display curvature when plotted against $1/n$.
We show the ratios plotted against $1/n$ and $1/\sqrt{n}$ in Figure~\ref{fig:r12}. 
 
\pagebreak
\clearpage

These plots are behaving as expected, with the plot against $1/n$ displaying considerable curvature, while the plot against $1/\sqrt{n}$ is
visually linear. This is strong evidence that $\beta=1/2$, just as  is the case for stack polyominoes. 

\begin{figure}[h!] 
\begin{minipage}[t]{0.43\textwidth} 
\centerline{\includegraphics[width=\textwidth]{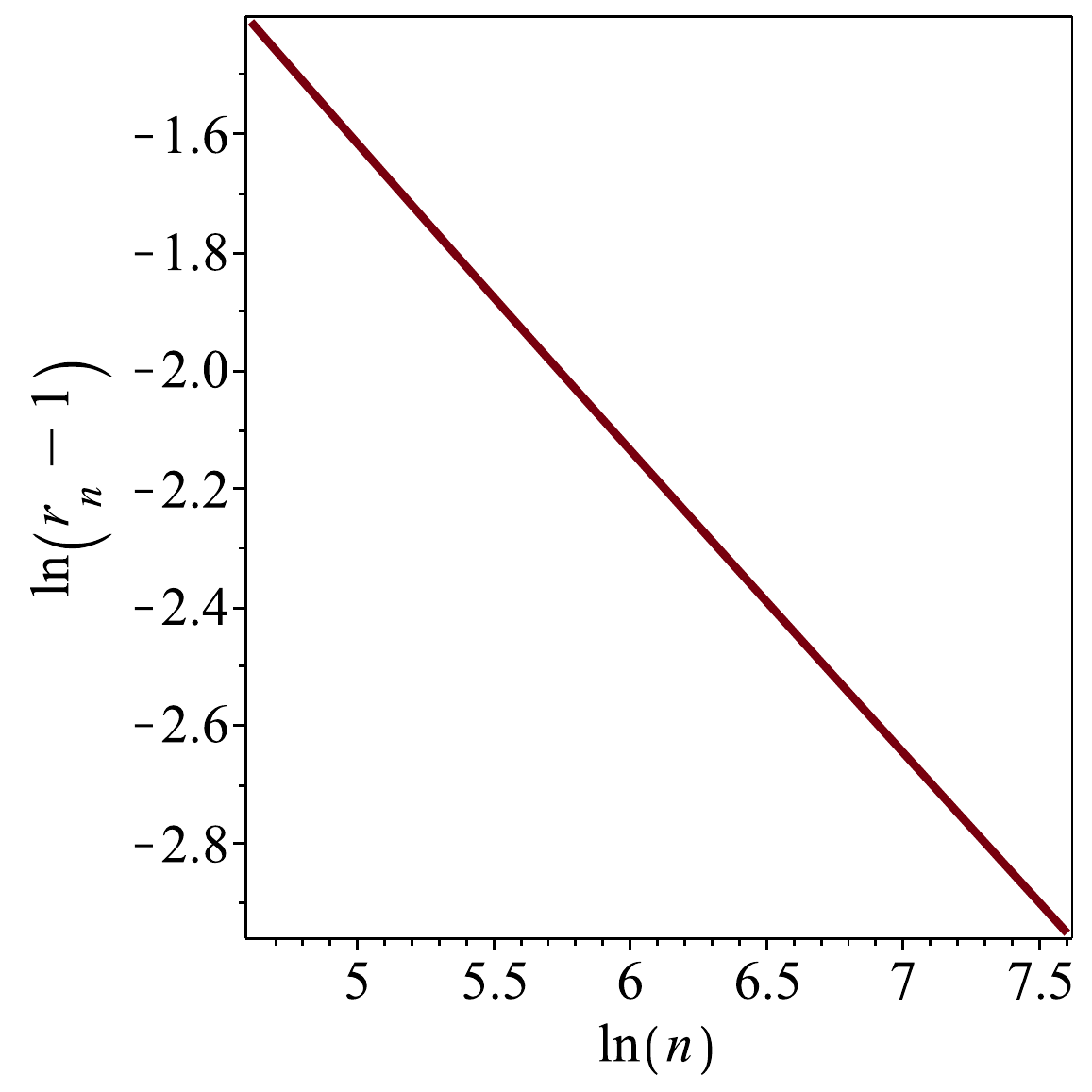} }
\caption{Log-log plot of $r_n-1$ against $n$.}
\label{fig:ll1}
\end{minipage}
\hspace{0.05\textwidth}
\begin{minipage}[t]{0.43\textwidth} 
\centerline{\includegraphics[width=\textwidth]{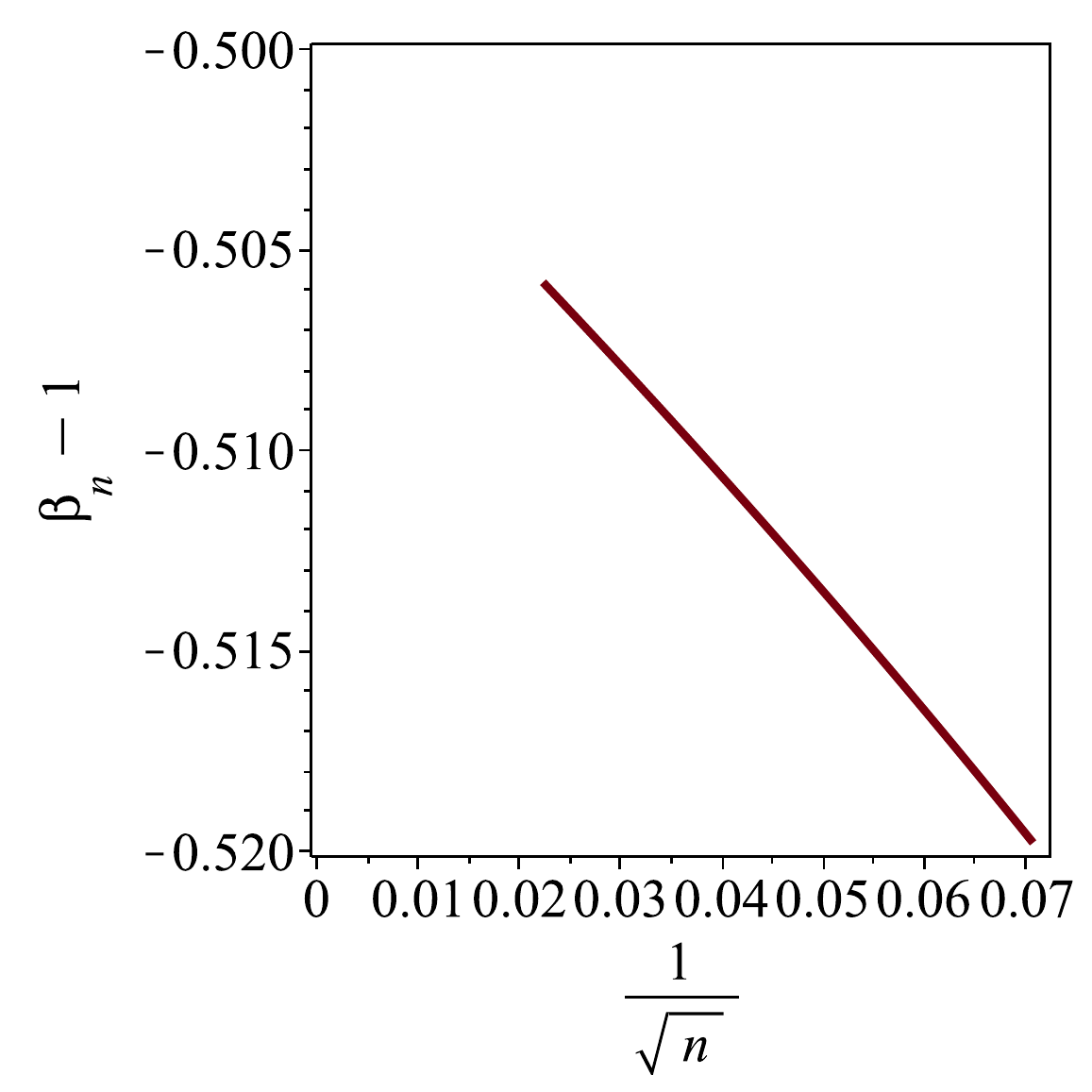}}
\caption{Gradient of log-log plot.}
\label{fig:grad}
\end{minipage}
\end{figure}

In fact we can easily refine this estimate.
From Equation~\eqref{eq:rn}, one sees that 
\BE \label{eq:rbeta}
r_n-1 = a\beta \pi \cdot n^{\beta-1} + O\left ( \frac{1}{n} \right ).
\EE
 Accordingly, a plot of $\log(r_n-1)$ versus $\log{n}$ should be linear, with gradient $\beta-1$. We would expect an estimate of $\beta$ close to that which linearised the ratio plot. In Figure~\ref{fig:ll1} we show the log-log plot, and in Figure~\ref{fig:grad} we show the local gradient plotted against $1/\sqrt{n}$. The linearity of the first plot is obvious, while the second is convincingly going to a limit of $-0.5$ as $n \to \infty$.

\begin{figure}[h!] 
\begin{minipage}[t]{0.49\textwidth} 
\centerline{\includegraphics[width=.85\textwidth,angle=0]{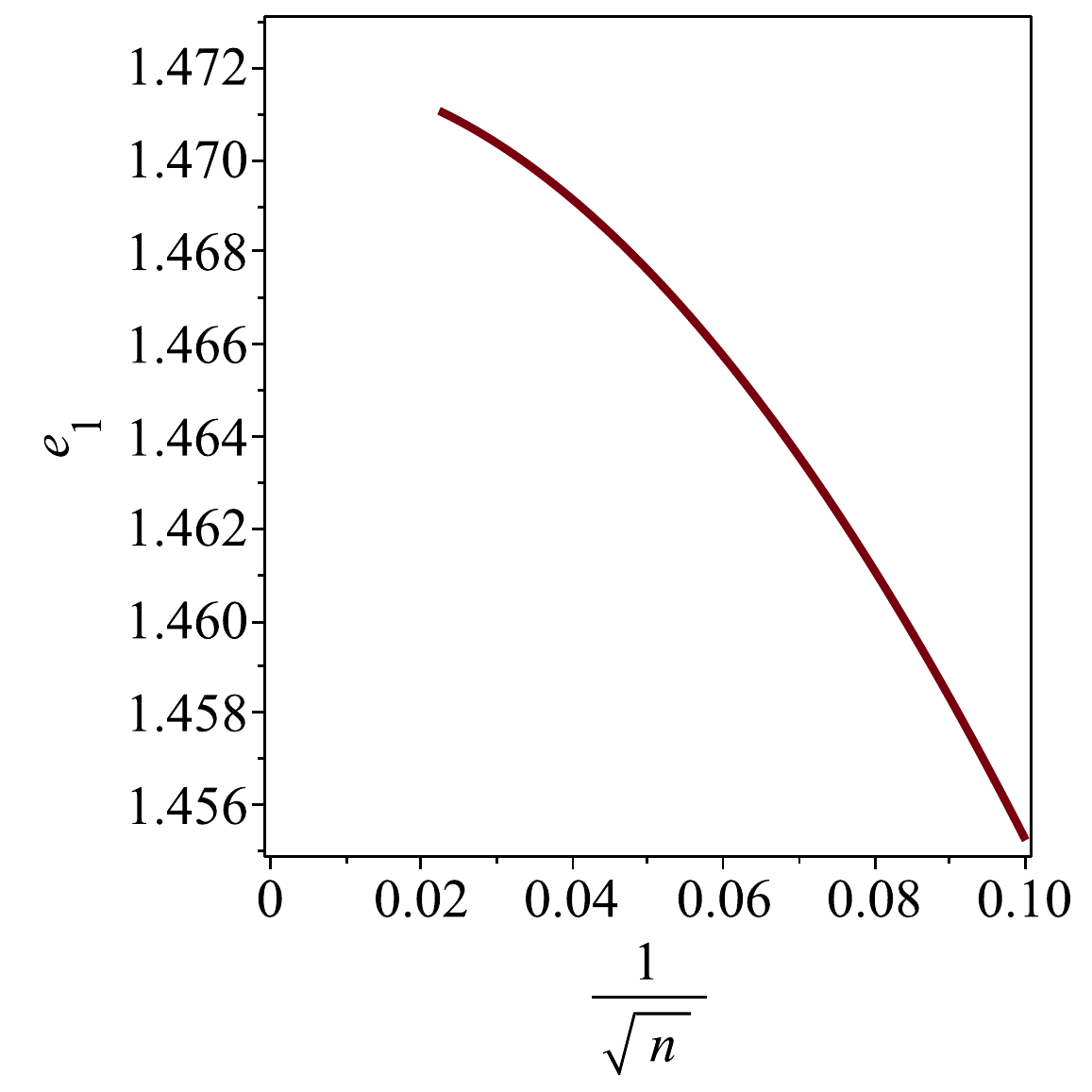} }
\setlength{\abovecaptionskip}{5pt} 
\caption{Plot of $e_1$ against $1/\sqrt{n}$.}
\label{fig:e1}
\end{minipage}
\begin{minipage}[t]{0.49\textwidth} 
\centerline{\includegraphics[width=.85\textwidth,angle=0]{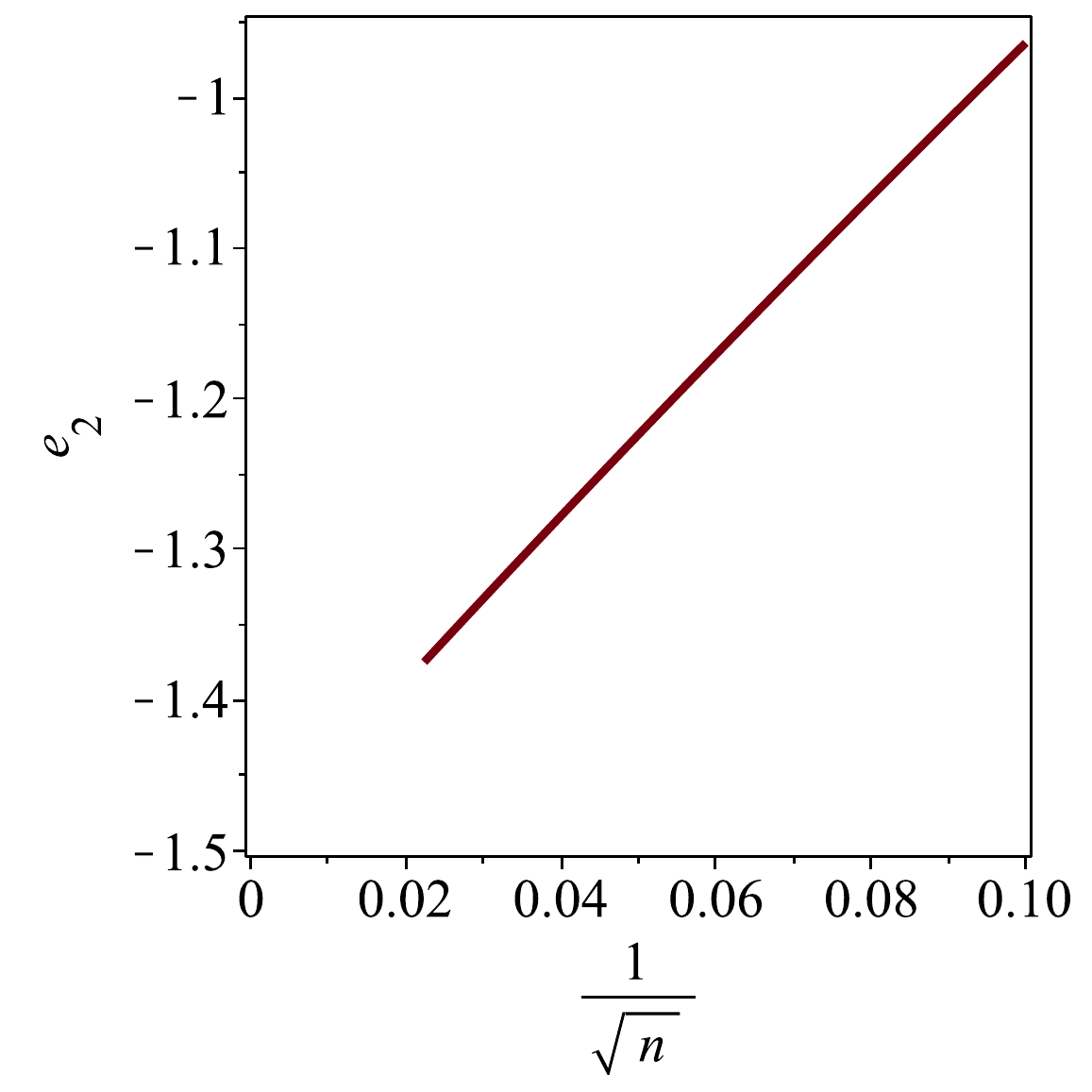}}
\setlength{\abovecaptionskip}{5pt} 
\caption{Plot of $e_2$ against $1/\sqrt{n}$.}
\label{fig:e2}
\end{minipage}
\end{figure}

Having convincingly established that $\beta=1/2,$ just as for stack polyominoes, it remains to determine the other parameters. There are several ways one might proceed, but here is one that works quite well. From the conjectured asymptotic form, we write
\BE \label{eqn:fit}
\lambda_n := \frac{\log(l_n)}{\pi\sqrt{n}}  \sim a -\frac{\delta\log{n}}{\pi\sqrt{n}} -\frac{\log{c}}{\pi\sqrt{n}} ,
\EE
so one can readily fit successive triple of coefficients $\lambda_{k-1}, \lambda_k, \lambda_{k+1},$ to the linear equation
\BE \lambda_n=e_1+e_2\frac{\log{n}}{\pi\sqrt{n}} +e_3\frac{1}{\pi\sqrt{n}}, \EE  with $k$ increasing until one runs out of known coefficients. Then $e_1$ should give an estimator of $a,$ $e_2$ should give an estimator of $-\delta$ and $e_3$ should give an estimator of $-\log(c)$. The result of doing this is shown for $e_1$ and $e_2$ in Figures~\ref{fig:e1} and~\ref{fig:e2} respectively. 

\pagebreak

We estimate the limits as $n \to \infty$ of $e_1$ as approximately $1.472,$ and $e_2$ as $-1.5$. From the asymptotic expression for $s_n,$ we expect $a$ to likely involve a square root. So we look at $e_1^2 = 2.16678,$ which we conjecture to be $13/6$. The exponent $\delta$ is expected to be a simple rational, and $3/2$ is indeed a simple rational! We don't show the plot for $e_3,$ as it does not give a precise enough estimate to conjecture the value of $\log(c)$ with any precision.

So at this stage we can reasonably conjecture that \BE l_n  \sim \frac{\exp(\pi\sqrt{13n/6})}{c \cdot n^{3/2}}. \label{conj1}\EE  We reached this stage based on only 100 terms in the expansion. In order to both gain more confidence in the conjectured form, and to calculate the constant, we needed more terms, and eventually generated 2000 terms from Equation~\eqref{eq:lconv}. 

With hindsight, an arguably more elegant way to analyse this series is to consider only the coefficients $l_{n^2}$. Denote $\ell_n:=l_{n^2}$. The conjectured form~\eqref{conj1} then  becomes  $\ell_n  \sim \frac{\exp(n\pi\sqrt{13/6})}{c \cdot n^{3}}$.  We have 44 coefficients of the series $\ell_n$ available, and these grow in the usual power-law manner, that is,
$\ell_n \sim D\cdot \mu^n \cdot n^g$.

\begin{figure}[h!] 
\begin{minipage}[t]{0.49\textwidth} 
\centerline{\includegraphics[width=\textwidth]{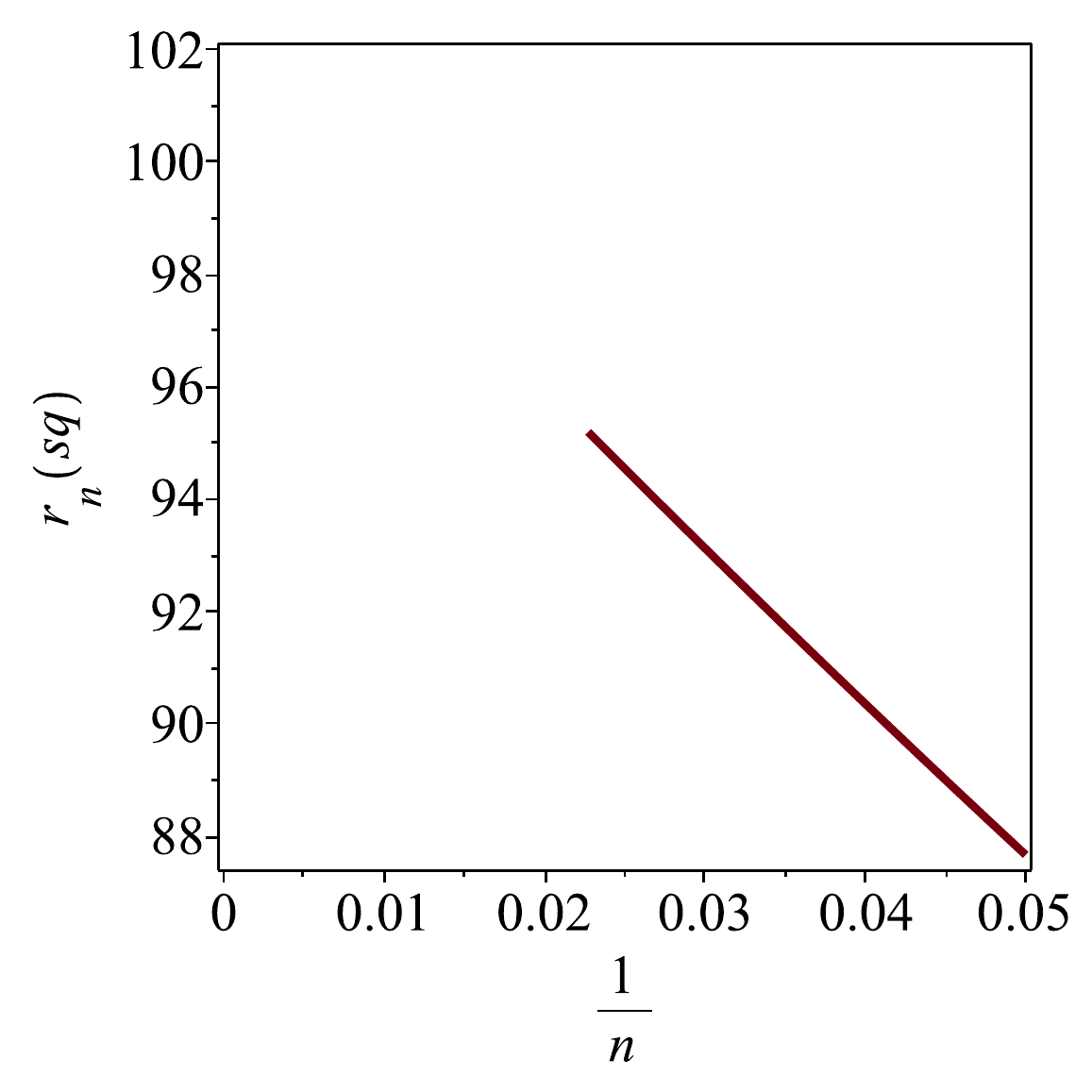} }
\setlength{\abovecaptionskip}{3pt} 
\caption{Plot of ratios $r_n^{(sq)}$ against $1/{n}$.}
\label{fig:r3}
\end{minipage}
\begin{minipage}[t]{0.49\textwidth} 
\centerline{\includegraphics[width=\textwidth]{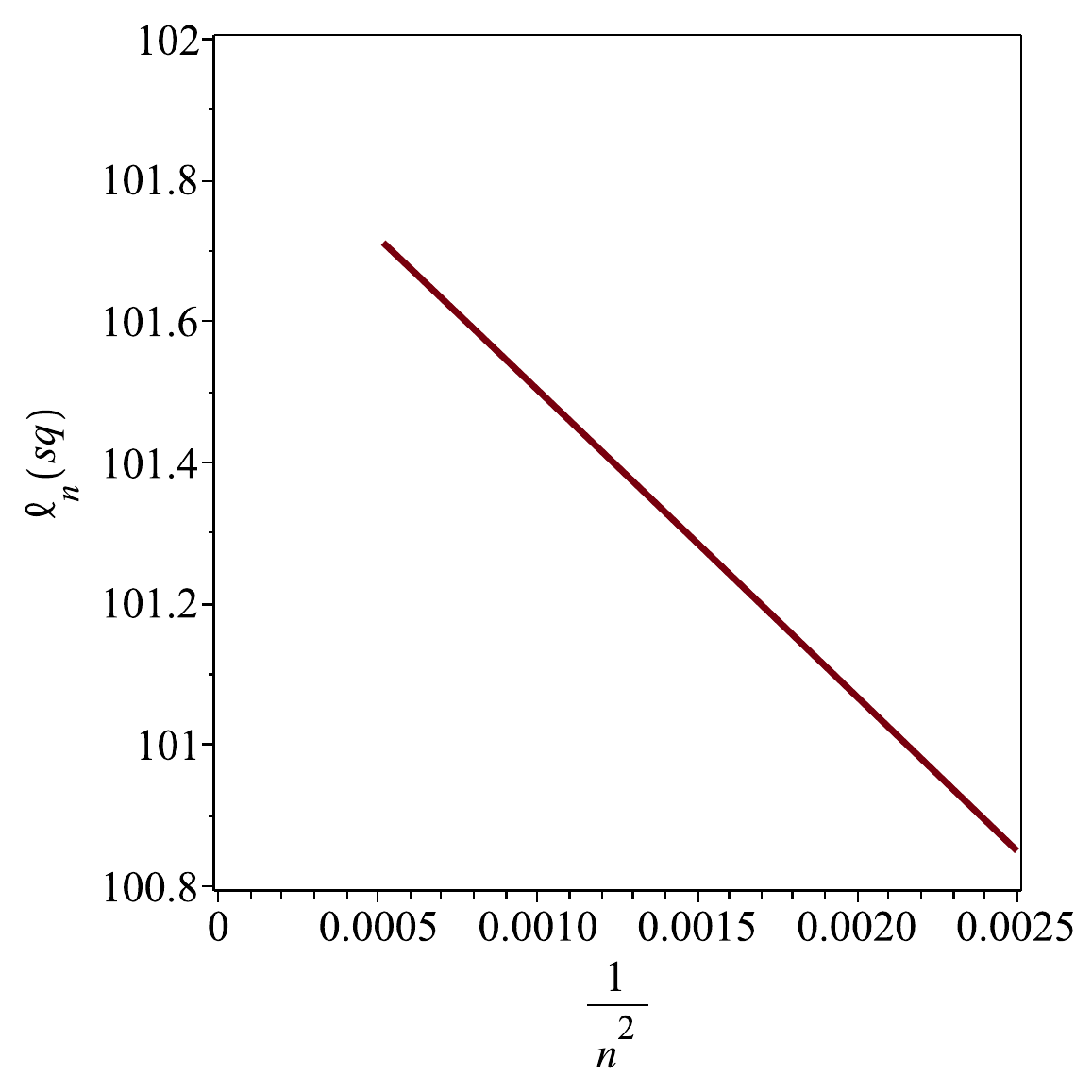}}
\setlength{\abovecaptionskip}{3pt} 
\caption{Plot of linear intercepts $\ell_n^{(sq)}$ against $1/{n^2}$.}
\label{fig:r4}
\end{minipage}
\end{figure}

We now analyse this sequence 
assuming its asymptotic form to be \BE \ell_n  \sim \frac{\exp(n\pi\sqrt{a})}{c \cdot n^{b}},\label{asymptform} \EE  with $a,$  $b,$ and $c$ to be determined.
Then we  form the ratios, \BE r_n^{(sq)}=\ell_n/\ell_{n-1} = \mu(1-b/n+ o(1/n)),\EE  where $\mu=\exp(\pi \sqrt{a})$ and where the superscript $(sq)$ is a mnemonic recalling that we here consider sequences derived from the subsequence $\ell_n=l_{n^2}$  of square indices. 
 Plotting the ratios $r_n^{(sq)}$ against $1/n$ should give a linear plot with gradient $-b\mu$ and ordinate $\mu$.  For a pure power-law the term $o(1/n)$ is $O(1/n^2),$ and the estimate of $\mu$ can thus be refined by plotting the linear intercepts $\ell_n^{(sq)}=n\cdot r_n-(n-1)\cdot r_{n-1}$ against $1/n^2$. The results of doing this are shown in Figures~\ref{fig:r3} and~\ref{fig:r4} for the ratios and linear intercepts respectively.
\pagebreak

\begin{figure}[h!]
\begin{minipage}[t]{0.49\textwidth} 
\centerline{\includegraphics[width=.9\textwidth]{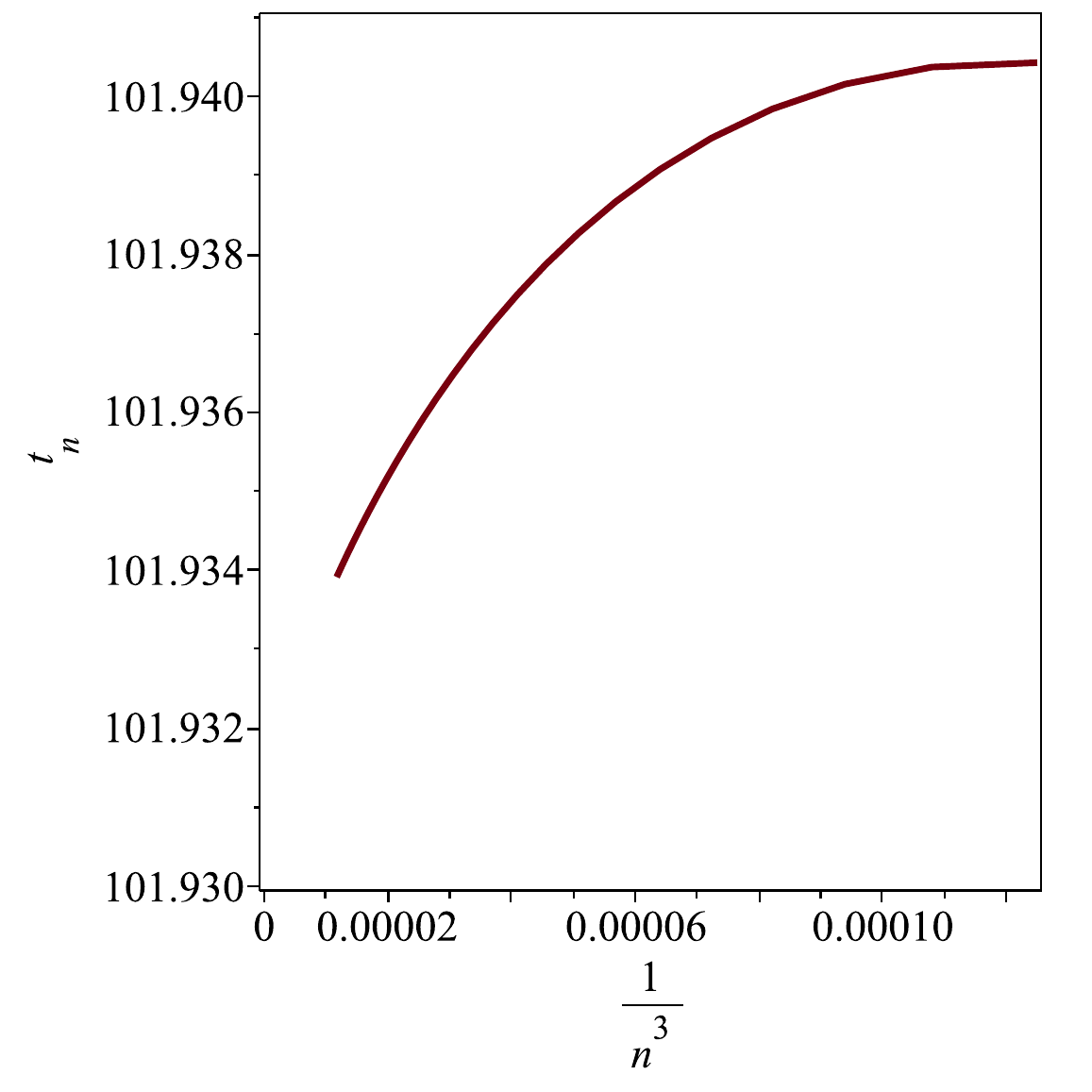} }
\setlength{\abovecaptionskip}{2pt} 
\caption{Plot of sequence $t_n$ against $1/{n^3}$.}
\label{fig:r5}
\end{minipage}
\begin{minipage}[t]{0.49\textwidth} 
\centerline{\includegraphics[width=.9\textwidth]{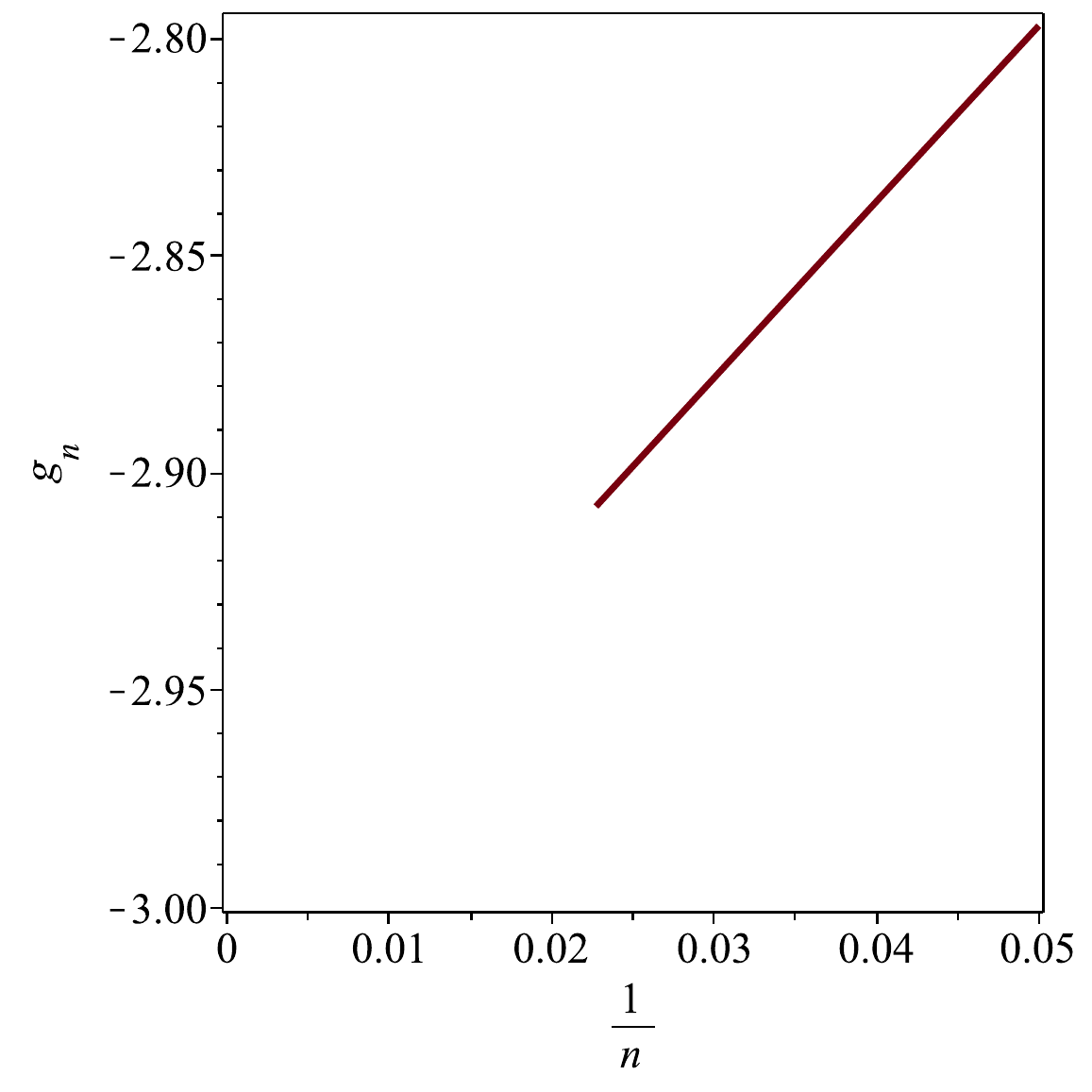}}
\setlength{\abovecaptionskip}{2pt} 
\caption{Plot of exponent estimates $g_n$ against $1/{n}$.}
\label{fig:g2}
\end{minipage}
\end{figure}

It can be seen that the linear intercepts have a faster convergence. We can go further and eliminate the $O(1/n^2)$ term by forming the sequence $t_n=(n^2\cdot \ell_n^{(sq)}-(n-1)^2 \cdot \ell_{n-1}^{(sq)})/(2n-1),$ and these are shown in Figure~\ref{fig:r5}. 
From this we estimate that the intercept of the plot with the ordinate is about 101.931. This is the growth constant $\mu=\exp(\pi\sqrt{a})$, from which we find $a \approx 2.16666,$ which strongly suggests that $a=13/6$ exactly. 

To estimate the exponent $b$ in the asymptotic form~\eqref{asymptform},
we introduce \BE g_n=(r_n^{(sq)}/\mu-1)\cdot n,\EE  noting that $\lim_n g_n= g =-b$.
Then, using the estimate of $\mu$ just given, we obtain the plot shown in Figure~\ref{fig:g2}. This is rather convincingly approaching $g=-3$.

\begin{figure}[h!] 
\centerline{\includegraphics[width=8.5cm,angle=0]{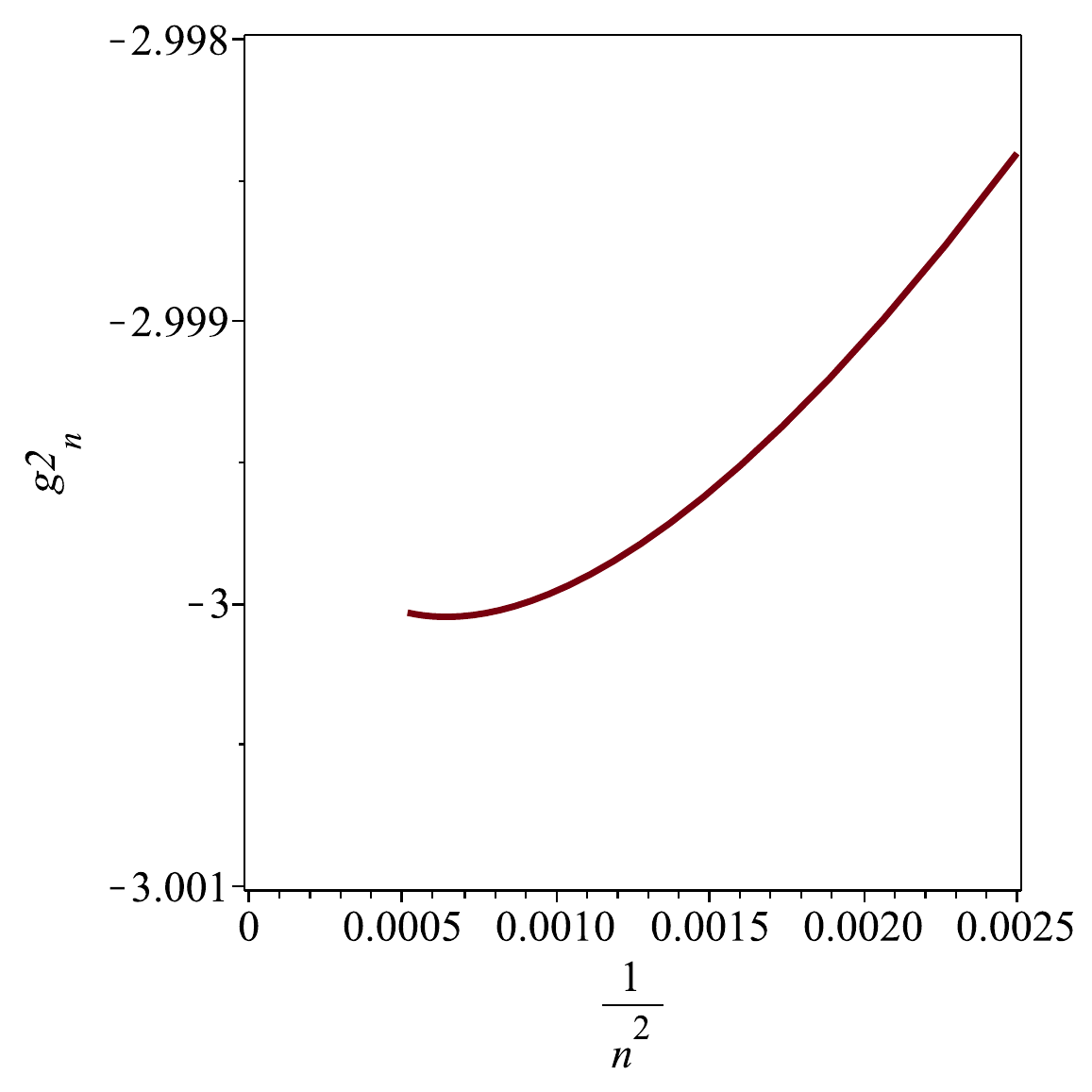} }
\setlength{\abovecaptionskip}{5pt} 
\caption{Plot of sequence $g2_n$ against $1/{n^3}$.}
\setlength{\abovecaptionskip}{5pt} 
\label{fig:g3}
\end{figure}

\clearpage
\noindent
We can do better by calculating the linear intercepts 
$g2_n:=n\cdot g_n - (n-1)\cdot g_{n-1}$. 
A plot of $g2_n$ against $1/n^2$ is shown in Figure~\ref{fig:g3}. The result $g=-3$ is totally convincing.

In order to calculate the constant $c$ in the asymptotic form~\eqref{asymptform}, we introduce the sequence \begin{equation}c_n:= \frac{\exp(\pi\sqrt{13n/6})}{l_n\cdot n^{3/2}},\end{equation} and extrapolate the sequence $c_n$ using any of a variety of standard methods. 

For this extrapolation, we used the Bulirsch--Stoer method (see \cite[Chapter~3.5]{SB80} or~\cite{BulirschStoer64} for more details), applied to the coefficient sequence $\{\ell_n\},$ with parameter $1/2,$ and $44$ terms in the sequence (corresponding to $44^2=1936$ terms in the original series). This gave the estimate $c \approx 0.023938510821419$. This unknown number is likely to involve a square root, cube root or fourth root of a small integer, just as did $s_n$. 

We investigated this by dividing by various powers of small integers, and tried to identify the result. Fortuitously, dividing the approximate value by $\sqrt{2}$ gave a result that the Maple command {\em identify} reported as $13/768$. This implies $c=13\sqrt{2}/768=0.023938510821419577\ldots,$ agreeing to all quoted digits with the approximate value. The occurrence of  $13$ in this fraction, as well as in the exponent square-root, is a reassuring feature, as is the factorisation of 768 as $3\cdot 2^8$.

Thus we conclude with the confident conjecture that the asymptotic form of the coefficients of $L$-convex polyominoes is
$$l_n  \sim \frac{3\cdot 2^8\cdot \exp(\pi\sqrt{13n/6})}{13\sqrt{2}\cdot n^{3/2}}.$$
\pagebreak

\section{201-avoiding ascent sequence}\label{201}
From the coefficients $u(n)$ for $n=0,\dots,27$ (this is the sequence~\oeis{A202062} in the OEIS),
we used the {\em gfun} package of Maple~\cite{SalvyZimmermann94}   and immediately found that the coefficients satisfy the recurrence relation
 \begin{align}
 \notag
 &(2n^2 + n)u(n) + (6n^2 + 45n + 60)u(n + 1) + (-34n^2 - 263n - 480)u(n + 2)\\ 
 \notag
 &+ (44n^2 + 421n + 984)u(n + 3) + (-20n^2 - 235n - 684)u(n + 4)\\
 \notag
 & + (2n^2 + 31n + 120)u(n + 5)=0, \\
 &{\rm with\,\,\,}u(0) = 1, u(1) = 1, u(2) = 2, u(3) = 5, u(4) = 15.
 \end{align}
 This recurrence can be converted to a second-order inhomogeneous ODE, 
or, as we prefer, a third-order homogeneous ODE, using the gfun command {\em diffeqtohomdiffeq}, giving
 \BE P_3(x)f^{'''}(x)+P_2(x)f^{''}(x)+P_1(x)f^{'}(x)+P_0(x)f(x)=0, \EE
where \\
\noindent\resizebox{\textwidth}{!}{\noindent
\begin{minipage}{1.144\textwidth}
\begin{align*}
 P_3(x) &=-2x^2(x^3 + 5x^2 - 8x + 1)(4x^4 - 30x^3 + 48x^2 - 36x + 15)(x - 1)^2,\\
 P_2(x)&=-3x(x - 1)(12x^8 - 30x^7 - 652x^6 + 2734x^5 - 4767x^4 + 4758x^3 - 2843x^2 + 870x - 85),\\
 P_1(x)&=-24x^9 + 30x^8 + 2754x^7 - 13278x^6 + 28884x^5 - 38106x^4 + 32436x^3 - 16620x^2 + 4350x - 420,\\
 P_0(x)&=\ 30(3x - 2)(3x^5 - 10x^4 + 19x^3 - 28x^2 + 24x - 7).
 \end{align*}
\end{minipage}} 
\smallskip

  The smallest root of the cubic factor in $P_3(x)$ is $x=0.1370633395\ldots$ and is the radius of convergence of the solution. Accordingly, the  growth constant $\mu$ thus satisfies
$$\mu=\frac{1}{x}=\frac{14}{3}\cos\left (  \frac{\arccos(\frac{13}{14})} {3} \right ) + \frac{8}{3}=7.295896946 \ldots$$
 
  This ODE can then be studied using the Maple  package {\em  DEtools}. We first convert the ODE to differential operator form through the command {\em de2diffop},   then factor this into the direct sum of two differential operators by the command {\em DFactorLCLM}. One of these operators is first order and one is second order. 
  
  The solution of the first order ODE is immediately given by the {\em dsolve} command, and is the rational function
 \begin{equation} \label{eqn:ODE1}
 y_1(x) = \frac{x^4 + 26x^3 - 45x^2 + 18x + 1}{12(x - 1)x^3}.
 \end{equation}

 To solve the second-order ODE, we obtain a series solution, the first term of which is O$(x^{-3})$. We multiply the solution by $x^3$ to obtain a regular power series, then use the {\em gfun} command {\em seriestoalgeq} to discover the cubic equation,
\begin{align} \label{cubic}
 4&(x-1)^3y_2(x)^3  \notag\\
 &-3(x - 1)(x^2 - x + 1)(x^6 - 235x^5 + 1430x^4 - 1695x^3 + 270x^2 + 229x + 1)y_2(x) \notag\\
 &+x^{12} + 510x^{11} - 14631x^{10} + 80090x^9 - 218058x^8 + 316290x^7 - 253239x^6 \notag\\
 & + 131562x^5 - 70998x^4 + 37950x^3 - 8955x^2 - 522x + 1=0.
 \end{align}
\pagebreak 

 This can be solved by Maple's, {\em solve} command, giving three solutions. Inspection of their expansion reveals the appropriate one, and simplifying this gives the following rather cumbersome solution:
 Let
 \begin{align*}
 P_1&=x^{12} + 510x^{11} - 14631x^{10} + 80090x^9 - 218058x^8 + 316290x^7 - 253239x^6+ 131562x^5 \\
  & \quad \quad  - 70998x^4 + 37950x^3- 8955x^2 - 522x + 1-24\sqrt{3x(x - 1)(x^3 + 5x^2 - 8x + 1)^7} ,\\
P_2 &= (x^2 - x + 1)(x - 1)^4(x^6 - 235x^5 + 1430x^4 - 1695x^3 + 270x^2 + 229x + 1),\\
 P_3 &= (3^{5/6}i + 3^{1/3}), \text{\ and\ } P_4=(-3^{5/6}i + 3^{1/3}). \text{ Then, one has}
 \end{align*}
 \begin{equation}
 y_2(x)=\frac{-3^{2/3}\left( P_4\left (-P_1\cdot (x-1)^6\right )^{2/3}+P_2\cdot P_3 \right)}{12\, \left(-P_1\cdot (x-1)^6\right)^{1/3}(x-1)^3}.
 \end{equation}
 The solution to the original ODE is then
 $$y(x)=\frac{y_2(x)}{12x^3}-y_1(x)=1+x+2x^2+5x^3+15x^4+\cdots$$
 
 This analysis required only the first 24 terms given in the OEIS, so the correct prediction of the next four terms gives us confidence that this is indeed the exact solution. 
 
 We next obtained the first 5000 terms in only a few minutes of computer time by expanding this solution.
 We used these terms to calculate the amplitude. That is to say, we now know that the coefficients behave asymptotically as
$u(n) \sim C \mu^n n^{-9/2}$.
 Equivalently, the generating function behaves as
 $$U(x)=\sum u(n)x^n =A(1-\mu \cdot x)^{7/2},$$
 where $C=A/\Gamma(-7/2) = 105A/(16\sqrt{\pi})$.
 We estimate $C$ by assuming a pure power law, so that $$\frac{u(n)\cdot n^{9/2}}{\mu^n} = C(1+\sum_{k \ge 1} a_k/n^k).$$ We calculated the first twenty coefficients of this expansion, which
 allowed us to estimate $C=13.4299960869\ldots$ with 74-digit accuracy (as checked later).
 Unless one is very fortunate (for example, when the Maple command {\em identify} determines an expression for this constant, which it doesn't in our case),  to identify this constant requires some experience-based guesswork. 
 
 Such constants in favourable cases are a product of rational numbers and square roots of small integers, sometimes with integer or half-integer powers of $\pi$.
 These powers of $\pi$ usually arise from the conversion factor in going from the generating function amplitude~$A$ to the coefficient amplitude~$C$. That is to say, we might expect the amplitude~$A$ to be simpler than $C$. And, to eliminate square-roots, we will try and identify $A^2$ rather than $A$.

 We do this by seeking the minimal polynomial with root $A^2$, using the command {\em MinimalPolynomial} in either Maple or Mathematica. 
In fact, one only requires 20 digit accuracy in the estimate of $A^2$ to
 establish the minimal polynomial, $A^6 - 1369A^4 + 17839A^2 + 1,$ which can be solved to give
 \begin{equation}C=\frac{35}{16}\left (\frac{4107}{\pi} - \frac{84}{\pi}\sqrt{9289}\cos\left (\frac{\pi}{3} +\frac{1}{3} \arccos\left [\frac{255709\sqrt{9289}}{24653006} \right ]\right )\right )^{1/2}.\end{equation}
\pagebreak
 
This derivation includes a degree of hindsight. In fact we searched for the minimal polynomial for the amplitude~$C$, by including various powers of $\pi,$ and then choose the polynomial of minimal degree. This required a much greater degree of precision in our estimate of $C$ to ensure we found the correct minimal polynomial.

It has been pointed out to us by Jean-Marie Maillard that the amplitude~$A$ can be obtained directly from the solution of the cubic equation~\eqref{cubic}, by extracting the coefficient of $(1-\mu \cdot x)^{7/2}$, as explained in~\cite[Chapter~VII.7.1]{FlSe09}. This gives the minimal polynomial that we obtained by numerical experimentation.  This alternative way 
to derive asymptotic expansions  is a more elegant method, as it is automatic, but it only works for algebraic functions.
There are thus many sequences for which it is not applicable,
as in the case of $L$-convex polyominoes (for which the generating function is not algebraic, as its radius of convergence is not algebraic), while our numerical approach can still yield conjecturally exact results. 
\vspace{-4mm}

\section{Conclusion}
We have shown how experimental mathematics can be used to conjecture exact asymptotics, in the case of $L$-convex polyominoes, and to conjecture an exact solution, in the case of $201$-avoiding ascent sequences. 
We hope that the results will be of interest, and that the methods will be more widely applied, as there are many outstanding combinatorial problems that lend themselves to such an approach. 

We recognise that these results are conjectural. 
We leave proofs to those more capable, and in the hope that the maxim of the late lamented J.~M.~Hammersley to the effect~that ``it is much easier to prove something when you know that it is true'' will aid that endeavour.

\subsection*{Acknowledgements.}
AJG wishes to acknowledge helpful discussions with Jean-Marie Maillard, and with Paolo Massazza on the topic of $L$-convex polyominoes, and to thank the\linebreak ARC Centre of Excellence for Mathematical and Statistical Frontiers (ACEMS) for support.
\vspace{-6.7mm}

\bibliographystyle{SLC} 
\def\UrlFont{\rmfamily}  
\bibliography{GuttmannKotesovec}
\end{document}